\documentclass[12pt,a4paper,leqno]{amsart}
\usepackage[utf8]{inputenc}
\usepackage{amssymb,color,tikz, cancel}
\usepackage{amsmath,amscd,amsfonts,amsthm,verbatim}
\usepackage{tikz-cd}
\usepackage{enumitem}
\usepackage[all]{xy}

\newcommand{\CC}{{\mathbb C}}
\newcommand{\QQ}{{\mathbb Q}}

\newcommand{\cE}{{\mathcal E}}
\newcommand{\cS}{{\mathcal S}}
\newcommand{\cF}{{\mathcal F}}
\newcommand{\cG}{{\mathcal G}}
\newcommand{\cP}{{\mathcal P}}

\newcommand{\cX}{{\mathcal X}}

\newcommand{\eps}{\varepsilon}
\newcommand{\cO}{{\mathcal O}}
\newcommand{\cM}{{\mathcal M}}

\newcommand {\PP}{\mathbb{P}}
\newcommand{\cJ}{{\mathcal J}}

\newcommand{\cL}{{\mathcal L}}
\newcommand {\EE}{\mathbb{E}}
\newcommand {\ZZ}{\mathbb{Z}}
\newcommand {\NN}{\mathbb{N}}

\newcommand{\isom}{\cong}

\DeclareMathOperator{\Gm}{\mathbb{G}_{\text{m}}}

\DeclareMathOperator{\EExt}{{\cE}xt}
\DeclareMathOperator{\Hilb}{Hilb}

\DeclareMathOperator{\Proj}{Proj}
\DeclareMathOperator{\Spl}{Spl}
\DeclareMathOperator{\PSpl}{PSpl}
\DeclareMathOperator{\PSSpl}{P{\cS}pl}
\DeclareMathOperator{\SSpl}{{\cS}pl}
\DeclareMathOperator{\Pic}{Pic}

\DeclareMathOperator{\End}{End}

\DeclareMathOperator{\tr}{tr}
\DeclareMathOperator{\ad}{ad}

\DeclareMathOperator{\rk}{rk}
\DeclareMathOperator{\rank}{rank}
\DeclareMathOperator{\id}{id}
\DeclareMathOperator{\Coh}{Coh}

\DeclareMathOperator{\Ext}{Ext}

\DeclareMathOperator{\Hom}{Hom}
\DeclareMathOperator{\Sing}{Sing}

\DeclareMathOperator{\Spec}{Spec}
\DeclareMathOperator{\Ker}{Ker}

\newtheorem{theorem}{Theorem}[section]
\newtheorem{lemma}[theorem]{Lemma}

\newtheorem{proposition}[theorem]{Proposition}
\newtheorem{corollary}[theorem]{Corollary}
 \theoremstyle{definition}
\newtheorem{definition}[theorem]{Definition}
\newtheorem{example}[theorem]{Example}
\newtheorem{remark}[theorem]{Remark}

\title[Simple sheaves on singular K3s]{On the moduli space of simple sheaves on singular K3 surfaces}

\date{\today}

\begin{document}

 \author[B.\ Fantechi]{Barbara Fantechi} 
 \address{SISSA
Via Bonomea 265, I-34136 Trieste, Italy}
  \email{fantechi@sissa.it, 
  ORCID 0000-0002-7109-6818}
  \author[R.\ M.\ Mir\'o-Roig]{Rosa M.\ Mir\'o-Roig} 
  \address{Facultat de
  Matem\`atiques i Inform\`atica, Universitat de Barcelona, Gran Via des les
  Corts Catalanes 585, 08007 Barcelona, Spain} \email{miro@ub.edu, ORCID 0000-0003-1375-6547}

\thanks{The first author has been partially supported by PRIN 2017 {\em Moduli theory and
birational classification} and PRIN 2022 {\em Geometry of algebraic structures: moduli, invariants, deformations}. She is a member of GNSAGA of INDAM} 
\thanks{The second author has been partially supported by the grant PID2020-113674GB-I00}

\begin{abstract} 
Mukai proved that the moduli space of simple sheaves on a smooth projective K3 surface is symplectic, and in \cite{FM2} we gave two constructions allowing one to construct new locally closed Lagrangian/isotropic subspaces of the moduli from old ones. 
In this paper, we extend both Mukai's result and our construction to reduced projective K3 surfaces; for the former we need to restrict our attention to perfect sheaves. There are two key points where we cannot get a straightforward generalization. In each, we need to prove that a certain differential form on the moduli space of simple, perfect sheaves vanishes, and we introduce a smoothability condition to complete the proof.
\end{abstract}

\maketitle

\tableofcontents

\section{Introduction}

In this paper, a K3 surface  will be a projective, Gorenstein, reduced and connected complex surface $X$ with $H^1(\cO_X)=0$ and $\omega_X\cong \cO_X$ trivial (see remark \ref{rem1} and section \ref{criteriasmoothability} for details and motivation of this definition).

Let $X$ be a  K3 surface, $L\in \Pic X$  and $r,c_2$ integers with $r\ge 1$; we denote by $\Spl(r;L,c_2)$ the (algebraic) moduli space of simple  sheaves on $X$ with determinant $L$, rank $r$ and second Chern class $c_2$.

The first result is an extension of Mukai's theorem that the moduli space of simple sheaves is symplectic, with symplectic form at a point $[F]$ given by Serre duality on the tangent space $\Ext^1(F,F)$. For this form of Serre duality to work, we need that $F$ is perfect, i.e. admits a finite locally free resolution. We denote by $\PSpl(r;L,c_2)$ the open subspace of $\Spl(r;L,c_2)$ parametrizing perfect sheaves. Note that every coherent sheaf on a smooth variety is perfect.

\begin{theorem}\label{Symplectic} Let $X$ be a  K3 surface, $L\in \Pic(X)$  and $r,c_2$ integers with $r\ge 1$.  
\begin{enumerate}
    \item The moduli space $\PSpl(r;L,c_2)$ is a smooth algebraic space of dimension $2rc_2-(r-1)L^2-2r^2+2$ (Proposition \ref{smoothness});
    \item Serre duality induces  a natural, nowhere degenerate $2$-form $\alpha $ on $\PSpl(r;L,c_2)$ (Proposition \ref{nondegenerate});
    \item If the pair $(X,L)$ is smoothable, then the form $\alpha $ is closed and $\PSpl(r;L,c_2)$ is a symplectic algebraic space (Theorem \ref{is_symplectic}).
\end{enumerate}
\end{theorem}

We remark that the notion of Lagrangian (locally closed) subspace can be defined for any space with a choice of nowhere degenerate $2$-form, whether it is closed or not. Therefore, we can speak about Lagrangian subspaces of $\PSpl(r;L,c_2)$.

Under suitable assumptions, in \cite{FM2} we give two methods to associate to a subspace $\cL$ of a moduli space of simple vector bundles on a smooth projective K3 surface another subspace $\cL'$ in a different moduli space of simple vector bundles; the first one is generalized  syzygy bundles, the second is extension bundles. The main result of that paper is that, for both constructions, $\cL$ is Lagrangian/isotropic if and only if $\cL'$ is.

As a first step towards generalizing these results to non smooth K3 surfaces, in section \ref{sec_lagrangian}, we show that the extension construction works for non smooth K3 surfaces and we note that generalized syzygy bundles were already defined in \cite{FM1} for any projective variety $X$ of dimension $\ge 2$, with no assumptions on the singularities; so no generalization is necessary.

The proof in \cite{FM2} that the extension construction sends Lagrangian (resp. isotropic) subspaces to Lagrangian (resp. isotropic) subspaces (and conversely) applies unchanged without the assumption that the K3 surface $X$ be smooth.

Let $v\ge 1$ be an integer. Denote by $U_e\subset \PSpl(r;L,c_2)$ the open locus parametrizing perfect simple vector bundles $F$ such that $H^0(F)=H^0(F^*)=0$ 
and assume $U_e\ne \emptyset $. We define 
$\cP_{U_e}=\{(F,V) \mid F\in U_e, V^*\subset H^1(F^*), \dim V=v\}$ and we denote by 
$\pi:\cP_{U_e}\to U_e$ the natural projection. For any smooth, irreducible, locally closed subspace  
$\cL\subset \PSpl(r;L,c_2)$, we set $\cP_{\cL}:=\pi ^{-1}(\cL\cap U_e)$ and let $\widetilde{ \cL}$ be the image of the locally closed embedding 
$\beta: \cP_{\cL}\to \PSpl(r+v;L,c_2)$  
induced by the set map 
$(F,V)\to E$ where $E$ is the extension bundle given by the exact sequence $0\to V\otimes \cO_X \to E \to F \to 0$. We have:

\begin{theorem} \label{lagrangian_extensions}
(1) The subspace $\cL$ is half dimensional in  $\PSpl(r;L,c_2)$ if and only if $\widetilde{\cL}$ is half dimensional in  $\PSpl(r+v;L,c_2)$;

 (2) The subspace $\widetilde{\cL}$ is isotropic (respectively, Lagrangian) if and only if $\cL$ is isotropic (respectively, Lagrangian).
\end{theorem}

The corresponding proof for the generalized syzygy bundles construction uses Dolbeault cohomology. Therefore, it does not work for non smooth K3 surfaces $X$ and we introduce a smoothability condition to
complete the proof.

Let $w\ge r+2$ be an integer. Denote by $U\subset \PSpl(r;L,c_2)$ the open locus parametrizing perfect simple vector bundles $F$ such that $H^1(F)=0$ 
and $F$ is globally generated. Assume $U\ne \emptyset $, define 
$\cG^0_{U}=\{(F,W) \mid F\in U, W\subset H^0(F), \dim W=w, W\otimes \cO_X\twoheadrightarrow F\}$ and 
let $\pi:\cG_U^0\to U $ be the natural projection. For any smooth, irreducible, locally closed subspace  
$\cL\subset \PSpl(r;L,c_2)$, we set $\cG^0_{\cL}:=\pi ^{-1}(\cL\cap U)$ and let $\cL'$ be the image of the locally closed embedding 
$\cG^0_{\cL}\to \PSpl(w-r;-L,L^2-c_2)$  
set-theoretically given by 
$(F,W)\to S:=\Ker(W\otimes \cO_X\to F)$. It holds:

\begin{theorem}\label{lagrangian} (1) The subspace $\cL$ is half dimensional in  $\PSpl(r;L,c_2)$ if and only if $\cL'$ is half dimensional in  $\PSpl(w-r;-L,L^2-c_2)$;

 (2) Assume that the pair $(X,L)$ is smoothable. Then $\cL'$ is isotropic (respectively, Lagrangian) if and only if $\cL$ is isotropic (respectively, Lagrangian).

\end{theorem}

Finally, in Section \ref{criteriasmoothability}, we discuss the smoothability condition for pairs $(X,L)$ where $X$ is a K3 surface and $L$ is an ample line bundle. We list examples of cases in which it is satisfied, and we give a sufficient criterion in case $X$ has isolated lci singularities and $L$ is ample.

\vskip 4mm
\noindent {\bf Acknowledgement.}  The first author is thankful for hospitality at Universitat de Barcelona and Henri Poincar\'e Institute, as well as to Sissa for a sabbatical year, during which part of this research took place. The first author would like to thank the participants at the conference in Nantes 2023, in particular Kieran O'Grady, for useful conversations.

\section{Notation and background material}

In this section, we gather some basic results that we will need later on and fix relevant notation. Throughout the paper we work over the complex numbers $\CC$. A {\em projective variety} $X$ is a connected, reduced, projective scheme of pure dimension $n$ over $\CC$.

\subsection{Smoothable  K3 surfaces and pairs}\label{smoothable_K3}
\begin{definition}\label{def_smoothable}    
    A (projective) {\em smoothing} of a projective variety $X$ is a flat projective morphism  $\pi:\cX\to B$ of relative dimension $n$ with $B$  a smooth connected curve together with  an isomorphism $\alpha:\cX_{b_0}\to X$ for a fixed   closed point $b_0\in B$ and such that for every $b\in B\setminus\{b_0\}$ the fiber $\cX_b$ is smooth.
    
    For $L\in Pic(X)$, a (projective) {\em smoothing of the pair} $(X,L)$ is a (projective) smoothing  $\pi:\cX\to B$ of $X$ plus a line bundle $\cL\ \in \Pic(\cX)$ such that $\cL_{|\cX_{b_0}}\cong \alpha ^*L.$ We will say that a pair $(X,L)$ is {\em smoothable} if it admits a (projective) smoothing.
\end{definition}

\begin{remark}
    Note that we use projective morphism in the sense of Hartshorne, that is we assume that there exists a relative very ample line bundle.
\end{remark}

By $X$ is a {\em K3 surface} we mean a projective, Gorenstein, reduced
and connected surface $X$ with $H^1(\cO_X ) = 0$ and $\omega_X \cong \cO_X.$

The following remark justify our definition of K3 surface.

\begin{remark} \label{rem1} (1) If a  projective, Gorenstein, reduced
and connected  surface $X$ has a (projective) smoothing with general fiber a smooth K3 surface then $\chi (\cO_X)=2$ and $\omega _X\cong \cO_X$ (because $\Pic(X)$ is separated). Therefore, $h^0(\cO_X)=h^2(\cO_X)=1$ and $H^1(\cO_X)=0$.

(2) Conversely, if a singular K3 surface $X$ admits a (projective) smoothing $\pi:\cX\to B$, then for generic $b\in B\setminus\{b_0\}$ one has $H^1(\cO_{\cX _b})=0$ (by semicontinuity) and $\omega _{\cX_b}\cong \cO_{\cX_b}$ (because $\Pic(\cX/B)\longrightarrow B$ is unramified). Therefore, for generic $b\in B\setminus\{b_0\}$, the fiber $\cX_b$ is a smooth K3 surface.
    \end{remark}

As an example of a smoothable pair $(X,L)$ with $X$ a K3 surface we have $(X,\cO_X(d))$ for any $d\in \ZZ$, where $X$ is
\begin{enumerate}
    \item a quartic surface in $\PP^3$;
    \item a complete intersection of type (2,3) in $\PP^4$;
    \item a complete intersection of type (2,2,2) in $\PP^5$;
    \item a double cover of $\PP^2$ branched over a sextic.
\end{enumerate}

More generally, examples can be drawn from Reid's list of complete intersection K3 surfaces in weighted projective spaces (see \cite[Proposition 2.1]{Y}). 

A full classification of such pairs is not known. In the last section of this paper, we give a sufficient condition for a pair $(X,H)$, with $X$ a K3 surface and $H$ an ample line bundle, to be smoothable.

\subsection{Moduli of simple perfect  sheaves on surfaces}

Let $X$ be a projective variety. A coherent sheaf $\cF$ on $X$ is called {\em perfect} if it is perfect as an element of $D^b(\Coh(X))$; equivalently, if it admits a finite resolution by locally free sheaves of finite rank. 

\begin{remark}
   For any flat family of sheaves over $X$, the locus of points corresponding to perfect sheaves is open; in other words, perfect sheaves define an open substack of the stack of coherent sheaves. 
\end{remark}

A perfect sheaf $\cF$ on $X$ has well defined rank, determinant line bundle and Chern classes $c_i(\cF)\in H^{2i}(X^{an},\mathbb Z)$ induced by those of vector bundles by taking any locally free resolution. Note that rank and Chern classes are constant for the fibers of a flat family  of perfect sheaves with connected basis $B$. The determinant defines a morphism $B\to \Pic(X)$ and thus is constant if $H^1(X,\cO_X)=0$, i.e., if $\Pic(X)$ is discrete.

A coherent sheaf $\cF$ on $X$ is called {\em simple} if  the natural injection $$H^0(X,\cO_{X}) \longrightarrow \End_{\cO_X}(\cF)$$ is an isomorphism.
\begin{proposition}
 Let $X$ be a K3 surface, $L\in Pic(X)$, and $r,c_2$ integers with $r\ge 1$. The set of simple sheaves on $X$ which are perfect of rank $r$, determinant $L$, and second Chern class $c_2$ defines an open substack in $\SSpl(r;c_1(L),c_2)$. 
\end{proposition}

\begin{proof}
    We already know that simple sheaves of fixed numerical first Chern class are an algebraic stack. Fixing the determinant is an open condition, because $H^1(X,\cO_X)=0$ and being perfect is also an open condition on flat families of sheaves.
\end{proof}

\begin{definition}
    \label{modSimplePerfect}

Let $X$ be a K3 surface; we can give the set of isomorphism classes of perfect simple sheaves of rank $r$, determinant $L$ and second Chern class $c_2$ a structure of algebraic space as follows.

We first define a stack $\PSSpl_X(r;L,c_2)$ by associating to each scheme $B$ the pairs $(\cF,\alpha)$ where $\cF\in \Coh(X\times B)$ is a flat family of perfect simple sheaves of given rank and $c_2$, and $\alpha:\det \cF\to p_X^*L$ is an isomorphism. An isomorphism between pairs $(\cF,\alpha)$ and $(\cG,\beta)$ is an isomorphism $\phi:\cF\to \cG\in \Coh(X\times B)$ such that $\alpha=\beta\circ \det\phi$. 
In particular, the automorphism group of $(\cF,\alpha)$ is naturally isomorphic to $\mu_r$ if $B$ is connected.

The stack $\PSSpl_X(r;L,c_2)$ is a DM algebraic stack, and we define $\PSpl_X(r;L,c_2)$ to be the corresponding (coarse moduli) algebraic space. 
We will omit the subscript $X$ when no confusion is likely to arise.
\end{definition}

\begin{remark}
(1) The structure morphism $$\PSSpl(r;L,c_2)\to \PSpl(r;L,c_2)$$ is a $\mu_r$ gerbe and the morphism $\PSSpl(r;L,c_2)\to \SSpl(r;L,c_2)$ given by forgetting $\alpha$ is a principal $\Gm$ bundle over its image, the open locus of perfect sheaves.

(2) For smooth K3 surfaces we have $\PSpl(r;L,c_2)=\Spl(r;c_1(L),c_2)$ while $\PSSpl(r;L,c_2)$ is a principal $\Gm$ bundle over $\SSpl(r;c_1(L),c_2)$.

(3) Simple vector bundles are by definition perfect.

\end{remark}

\subsection{Moduli of simple perfect sheaves in families}

    In this section we will fix $\pi:\cX \to B$ a flat projective family of K3 surfaces and a line bundle $\cL\in \Pic(\cX)$. We define a DM algebraic stack $\PSSpl_\pi(r;\cL,c_2)$ over $B$ by associating to every morphism of schemes $B_1\to B$ 
 the pairs $(\cF,\alpha)$ where $\cF\in \Coh(\cX_1)$ is a flat family of simple perfect sheaves with given rank and $c_2$, and $\alpha:\det \cF\to g^*\cL$ is an isomorphism; where $\cX_1$ and $g$ are defined by the  cartesian diagram \[
\xymatrix{ \cX_1\ar[r]^g \ar[d]^{\pi_1} & \cX \ar[d]^{\pi}\\
 B_1\ar[r]_f &B.
}
 \]

\begin{proposition}
    The algebraic stack $\PSSpl_\pi(r;\cL,c_2)$ pulls back under base change in $B$, in the sense that for every cartesian diagram as above
we have a natural isomorphism \[
    \PSSpl_{\pi_1}(r;g^*\cL,c_2)\cong \PSSpl_\pi(r;\cL,c_2)\times_BB_1.
    \]
\end{proposition}

This implies the same statement for the corresponding algebraic spaces $\PSpl_\pi(r;\cL,c_2)$. In particular, for every $b\in B$ we have a natural isomorphism \[
    \PSSpl_{X_b}(r;\cL|_{X_b},c_2)\cong \PSSpl_\pi(r;\cL,c_2)\times_B
    \{b\}.
    \]

\subsection{Ext groups and trace maps} 
We gather together the technical results that we need in next sections. All of them are well known for coherent sheaves on smooth projective varieties and we are going to extend them to perfect sheaves on non-necessarily smooth projective varieties.

\begin{lemma} Let $X$ be a projective variety and let $F$ be a perfect sheaf on $X$.

(1) There are natural maps $i:\cO_X\to {\mathcal E}nd(F)$ and $\tr:{\mathcal E}nd(F)\to \cO_X$ which induces morphisms $i:H^j(\cO_X)\to \Ext^j(F,F)$ and $\tr:\Ext^j(F,F)\to H^j(\cO_X)$ verifying:
\begin{itemize}
    \item[(a)] $\tr\circ i=\rank(F)\id_F$; and 
    \item[(b)] $\tr(\psi\circ \varphi)=(-1)^{j_1+j_2}\tr(\varphi \circ \psi)$  for any $\psi \in Ext^{j_1}(F,F)$ and $\varphi \in \Ext^{j_2}(F,F)$.
\end{itemize}

 (2) If $r>0$, then for every integer $i$, there is a canonical isomorphism $$\Ext^i(F,F)\cong \Ext^i(F,F)_0\oplus H^i(X,\cO_X)$$ where $\Ext^i(F,F)_0:=\Ker[tr:\Ext^i(F,F)\to H^i(\cO_X)]$. 
\end{lemma}

\begin{proof} (1) It follows step by step the proof of \cite[Lemma 10.1.3]{HL} where the hypothesis of $X$ being smooth is only used to assure that $F$ has a finite free resolution by locally free sheaves of finite rank which in our case is true by assumption.

(2) As in  \cite[Proposition 2]{Ar} it follows from the splitting of the trace map due to point (1)(a): $$0 \to \cO_X \rightleftarrows {\mathcal E}nd(F) \to {\mathcal E}nd(F)/\cO_X \to 0.$$ 
\end{proof}

We quickly sketch how to  generalize Serre's duality to perfect sheaves on non necessarily smooth Gorenstein projective varieties.
\begin{lemma} \label{SerreDuality} Let $X$ be a Gorenstein projective variety of dimension $n$. For perfect sheaves $F$ and $G$ on $X$ the natural pairing
$$
\Ext^i(F,G)\otimes \Ext^{n-i}(G,F\otimes \omega _X)\longrightarrow \Ext^n(F,F\otimes \omega _X)\xrightarrow{\,\,\,\tr\,\,\,} H^n(X,\omega _X)\cong \CC$$
is perfect. In particular, the induced natural map $$\alpha(G,F):\Ext^{n-i}(G,F\otimes \omega _X)\to \Ext^i(F,G)^*$$ is an isomorphism, functorial in both $F$ and $G$.
    \end{lemma}
\begin{proof} We proceed by induction on the projective dimension $pd(F)$ of $F$. If $pd(F)=0$, $F$ is locally free and the result is true. Assume $pd(F)>0$ and write 
 $$
0\longrightarrow F_1\longrightarrow E \longrightarrow F\longrightarrow  0$$
where $E$ is locally free and $pd(F_1)=pd(F)-1.$
The result follows from the following commutative diagram

{\small 
\[
\xymatrix{\Ext^{n-i}(G,F_1\otimes \omega _X)\ar[r] \ar[d]^{ \alpha (G,F_1) \,\, \,\ } &  \Ext^{n-i}(G,E\otimes \omega _X)\ar[r] \ar[d]^{ \alpha (G,E) \,\, \,\ } & \Ext^{n-i}(G,F\otimes \omega _X)\ar[r] \ar[d]^{ \alpha (G,F) \,\, \,\ } & \\
\Ext^{i}(F_1,G)^*\ar[r] & \Ext^{i}(E,G)^*\ar[r] & \Ext^{i}(F,G)^*\ar[r] & 
}
\]
}
{\small 
\[
\xymatrix{ \Ext^{n-i+1}(G,F_1\otimes \omega _X)\ar[r] \ar[d]^{ \alpha (G,F_1) \,\, \,\ } & \Ext^{n-i+1}(G,E\otimes \omega _X) \ar[d]^{ \alpha (G,E) \,\, \,\ }\\
 \Ext^{i-1}(F_1,G)^*\ar[r] & \Ext^{i-1}(E,G)^*
}
\]
}
\noindent with exact rows, commutative squares  and natural vertical maps. By the 5-lemma, since all vertical morphisms but $\alpha (G,F)$ are isomorphisms, so is $\alpha (G,F)$.
\end{proof}

When $X$ is a K3 surface we will write this duality without the sheaf $\omega_X$, assuming a trivialization has been chosen.

\vskip 4mm
For readers' convenience, we end this  subsection with a useful result on trace maps valid for arbitrary families of perfect complexes over flat projective families.

\begin{proposition}\label{rel_ad_and_tr}
    Given a flat projective morphism $p:Y\to S$ of schemes and a perfect sheaf $\cF\in \Coh(Y)$ flat over $S$, there are natural morphisms $c:\EExt^i_p(\cF,\cF)\otimes \EExt^j_p(\cF,\cF)\to \EExt^{i+j}_p(\cF,\cF)$, 
$\tr^i:\EExt_p^i(\cF,\cF)\to R^ip_*\cO _{Y}$ and 
$\ad:R^ip_*\cO_{Y}\to \EExt^i_p(\cF,\cF)$ satisfying $\tr \circ  \ad= \rk(\cF)\id$ and $\tr(\psi \circ \varphi)=(-1)^{\deg(\psi)\deg(\varphi)}\tr(\varphi\circ \psi)$.
\end{proposition}
\begin{proof}
    We can use the same argument as in \cite[Section 10.1.7]{HL}, since this only depends on having a finite free resolution. 
\end{proof}

\section{Symplectic structure on moduli of perfect simple sheaves on K3}

The goal of this section is to prove  Theorem \ref{Symplectic}: if $X$ is a K3 surface and $L\in \Pic(X)$, then the algebraic space $\PSpl(r;L,c_2)$ is smooth of dimension $2rc_2-(r-1)L^2-2r^2+2$ and it carries a nowhere degenerate 2-form. If, in addition, $(X,L)$ is smoothable, the 2-form is closed and  $\PSpl(r;L,c_2)$ is a symplectic algebraic space. As we will see the symplectic structure is induced pointwise by \[
\Ext^1(\cF,\cF)\otimes \Ext^1(\cF,\cF) \to \Ext^2(\cF,\cF)\to H^2(X,\cO_X)=\mathbb C.
\]
where the first map is the cup product, the second is induced by trace and is an isomorphism because $\cF$ is simple, and the final isomorphism is induced by a choice of a symplectic form on $X$, unique up to nonzero scalar.

\subsection{Moduli of perfect simple sheaf on singular K3 surfaces}
\begin{proposition}\label{smoothness} Let $X$ be a  K3 surface, $L\in Pic(X)$  and $r,c_2$ integers with $r\ge 1$.  
 Then, the moduli space $\PSpl(r;L,c_2)$ is a smooth algebraic space of dimension $2rc_2-(r-1)L^2-2r^2+2$.
    \end{proposition}
\begin{proof}

The tangent and obstruction spaces to $\PSpl(r;L,c_2)$ at a point $[F]$ are $\Ext^i(F,F)_0$ for $i=1,2$. Note that $\Ext^0(F,F)_0=0$ since $F$ is simple and $r\ge 1$.

By Lemma \ref{SerreDuality}, $\Ext^2(F,F)$ is dual to $\Ext^0(F,F)$ and it is one-dimensional. Therefore the relative obstruction $\Ext^2(F,F)_0$ is zero and the moduli space is smooth of dimension $\chi(\cO_X)-\chi(F,F)$; the explicit formula follows by Riemann-Roch, see \cite{Fu} Example 18.3.4(b).
\end{proof}

\begin{proposition}\label{nondegenerate}  Let $X$ be a  K3 surface, $L\in Pic(X)$  and $r,c_2$ integers with $r\ge 1$.  
 Then the moduli space  $\PSpl(r;L,c_2)$  has a natural nowhere degenerate 2-form $\alpha$ which at every point $F$ is the Serre duality  pairing\[
 \Ext^1(F,F)\times \Ext^1(F,F)\to H^2(\cO_X)\to \CC.
 \]
 \end{proposition}
\begin{proof}
The existence of a global $2$-form pointwise inducing the pairing in families follows from \cite[Proposition 10.3.5]{HL}; it is pointwise non-degenerate by Lemma \ref{SerreDuality}. The assumption that $X$ is smooth is only used to guarantee that $F$ has a finite locally free resolution.
\end{proof}

\begin{remark}
    Note that this result remains valid for any projective Gorenstein surface with $\omega_X\cong \cO_X$ because Serre duality induces a duality on $\Ext^1(F,F)_0$ (which is not equal to $\Ext^1(F,F)$ if $H^1(\cO_X)\ne 0$).
\end{remark}

To show that $\PSpl(r;L,c_2)$ is algebraic symplectic we still need to prove that $\alpha$ is closed. We want to reduce the proof of the closedness of $\alpha $ to the smooth case, assuming the pair $(X,L)$ is smoothable. In order to do so, we have to extend the construction of \cite[Proposition 10.3.5]{HL} to families of possibly singular surfaces.

\subsection{Symplectic structure on moduli of simple perfect sheaves on smoothable K3 surfaces}

 Let $\pi:\cX\to B$ and $\cL\in \Pic(\cX)$ be a smoothing of $(X,L)$. We consider the relative moduli stack $f:\cM=\PSSpl_\pi(r;\cL,c_2)\to B$.

\begin{lemma}
    The  morphism $f:\cM \to B$ is smooth of relative dimension $2rc_2-(r-1)L^2-2r^2+2$.
\end{lemma}
\begin{proof}
Let $p:=(b,F, \alpha)$ be a point in $\cM$, with $b\in B$, $F\in \Coh(\cX_b)$ and $\alpha:\det F\cong \cL|_{\cX_b}$ . The 
   relative tangent and obstruction spaces to the morphism $f$ at $p$ are $\Ext^1_{\cX_b}(F,F)_0$ and $\Ext^2_{\cX_b}(F,F)_0$; the result follows by the argument in Proposition \ref{smoothness}.
\end{proof}

We consider the cartesian diagram:
\[
\xymatrix{ \cX \times _B \cM\ar[r]^{g} \ar[d]^p & \cX \ar[d]^{\pi }\\
\cM\ar[r]^{f} & B.
}
\] Set  $\cX_{\cM}:=\cX \times _B \cM$ and let $\cF\in Coh(\cX_{\cM})$ be the universal sheaf.
Choose an \'etale  chart $S\to\cM$ with $S$ a scheme ($\cM$ has an \'etale atlas because it is  a Deligne-Mumford stack). Set $Y:=\cX_\cM\times_{\cM}S$, let $q:Y\to S$ be the natural projection and let $\cG$ be the pullback of $\cF$ to $Y$.

\begin{proposition}
    For $i=0,1,2$ the sheaves $\EExt^i_q(\cG,\cG)$ are locally free and commute with arbitrary base change. Moreover, $\EExt^1_q(\cG,\cG)$ is naturally isomorphic to the relative tangent bundle of $S$ over $B$.
\end{proposition}
\begin{proof}
  Let $s$ be any point in $S$ and $(b,F,\alpha)$ its image in $\cM$.   We apply \cite[Theorem 1.4]{L} first to $i=3$; and notice that $\Ext^3_{X_b}(F,F)=0$ because it is dual to $\Ext^{-1}_{X_b}(F,F)=0$. This implies that $\EExt^3_q(\cG,\cG)=0$, which is locally free. Hence, $\tau_2(s):\EExt^2_q(\cG,\cG)\otimes _S k(s)\to Ext^2_{X_b}(F,F)$ is surjective for any point $s\in S$. 
  
  We apply the theorem again for $i=2$, we get that for every $s$, the fiber of $\EExt^2_q(\cG,\cG)$ at $s$ is isomorphic to $\Ext^2_{X_b}(F,F)$ which is dual to  $\Hom_{X_b}(F,F)$ and hence is one dimensional. Therefore $\EExt^2_q(\cG,\cG)$ is locally free of rank $1$, since $S$ is smooth and thus reduced. By \cite[Theorem 1.4]{L}, $\tau_1$ is surjective at every $s\in S$.
  
  Again, the fiber of  $\EExt^1_q(\cG,\cG)$ at $s$ is isomorphic to $\Ext^1_{X_b}(F,F)$ which has dimension $2rc_2-(r-1)L^2-2r^2+2$. Since the dimension is the same at every $s\in S$,  $\EExt^1_q(\cG,\cG)$ is a vector bundle and commutes with base change.
  
  By the same argument also $\EExt^0_q(\cG,\cG)$ is a line bundle and commutes with base change.
  
  The last statement follows by $S\to \cM$ being \'etale, hence defining an isomorphism on tangent spaces.
  \end{proof}

\begin{remark}
       (1) By \'etale descent the locally free sheaves $\EExt^i_q(\cG,\cG)$ are pullbacks of similarly defined sheaves $\EExt^i_p(\cF,\cF)$ on $\cM$ and $\EExt^1_p(\cF,\cF)$ is naturally isomorphic to the relative tangent bundle $T_{\cM/B}$.

       (2) Let $\eps: \cM:=\PSSpl_\pi(r;\cL,c_2)\to \PSpl_\pi(r;\cL,c_2)$ be the structure map. Note that $f:\cM\to B$ factors uniquely via $\eps$ since $B$ is a scheme.
       
       Since each fiber $\Ext^i(F,F)$ of $\EExt^i_p(\cF,\cF)$ at a point $(b,F,\alpha)$ of $\cM$ is invariant under the $\mu_r$ action on $F$, it follows that the sheaves $\EExt^i_p(\cF,\cF)$ are canonicaly isomorphic to  $\eps^*\EE^i$, where $\EE^i:=\eps_*\EExt^i_p(\cF,\cF)$ are locally free on  $\PSpl_\pi(r;\cL,c_2)$. In particular, $\EE^1$ is isomorphic to  the relative tangent bundle of $\PSpl_\pi(r;\cL,c_2)$ over $B$ since $\eps$ is \'etale.
\end{remark}

\begin{theorem}
\label{global2form} The relative moduli space $\PSpl_\pi(r;\cL,c_2)$  has a natural  relative $2$-form $\alpha_B$ over $B$ whose restriction to each fiber $\PSpl(r;\cL|_{X_b},c_2)$ is the one constructed in Proposition \ref{nondegenerate}.
\end{theorem} 

\begin{proof}
    We first define a skew symmetric map $T_{f}\otimes T_f\to \cO_{\cM}$ by using the isomorphism $T_f\isom \EExt^1_p(\cF,\cF)$. The form $\alpha_B$ is obtained by composing the morphism $c:\EExt^1_p(\cF,\cF)\otimes \EExt^1_p(\cF,\cF)\to \EExt^2_p(\cF,\cF)$ with 
$\tr^2:\EExt_p^2(\cF,\cF)\to R^2p_*\cO_{\cX_\cM}$ and the isomorphism $R^2p_*\cO_{\cX_\cM}\isom \cO_{\cM}$ induced by a trivialization of $\omega_f\isom \cO_{\cX_\cM}$.
The restriction to each fiber follows from the naturality of the $c$ and $\tr$ constructions.
Finally, again everything descends to the coarse moduli space because the Ext groups are invariant under multiplication by $\mu_r$.
\end{proof}
\begin{theorem}\label{is_symplectic}
 In the same assumptions of Proposition \ref{nondegenerate}, assume moreover that the pair $(X,L)$ is smoothable. Then $\alpha$ is closed and makes $\PSpl(r;L,c_2)$  into a symplectic algebraic space.
\end{theorem} 
   
\begin{proof} 
Let $\pi:\cX\to B$ and $\cL\in \Pic(\cX)$ be a smoothing of $(X,L)$. We consider the relative moduli space $\PSpl_\pi(r;\cL,c_2)$, which is smooth over $B$ since the relative tangent space at each point $(b,F\in \PSpl(r;\cL|_{X_b},c_2))$ is again $\Ext^2(F,F)_0$ hence vanishes. We  construct on $\PSpl_\pi(r;\cL,c_2)$ a global $2$-form $\alpha_B$ (see Proposition \ref{global2form}) with the pointwise definition as in Proposition \ref{nondegenerate}, using that the relative tangent space at $(b,F)$ is $\Ext^1(F,F)_0$.

We consider the relative $3$-form $d_\pi\alpha_B$. It is zero over the locus where $X_b$ is smooth, which is dense in $B$. Therefore, it must be zero everywhere on $\PSpl_\pi(r;\cL,c_2)$.

Restricting $\alpha_B$ to the fiber at $b_0$ and using the isomorphisms $\cX_{b_0}\to X$ and $\cL|_{\cX_{b_0}}\to L$ concludes the argument.
\end{proof}

\section{Lagrangian submanifolds of  $\PSpl(r;L,c_2)$}\label{sec_lagrangian}

This section is entirely devoted to prove Theorem \ref{lagrangian_extensions} and Theorem \ref{lagrangian}. More precisely,  we show how using generalized syzygy bundles and extension bundles we can  associate to every Lagrangian
(resp. isotropic) algebraic subspace of $\PSpl(r; L,c_2)$ an induced Lagrangian (resp. isotropic) algebraic subspace of a different component of the moduli of perfect simple sheaves. 

\subsection{Extension bundles}

For sake of completeness we first quickly review from \cite{FM2} the construction of extension bundles.   
\begin{definition}\label{DefExt}
Let $F$ be a  vector bundle of rank $r$, determinant $L$ and second Chern class $c_2$ on a K3 surface $X$. For any subspace $V^*\subset  H^1(X,F^*)\cong \Ext^1(F,\cO_X)$ of dimension $v:=\dim V\ge 1$ we define the {\em extension  bundle} $E$ on $X$ associated to the couple $(F,V)$
as the vector bundle that comes up from the induced extension of $F$ by $V\otimes \cO_X$. Therefore, $ E$  fits into a short exact sequence:
\begin{equation}
    \label{extension1}
e: \quad 0\to V\otimes \cO_X\to E \to F\to 0
\end{equation}
where $e\in H^1(V\otimes F^*)=Hom(V^*,H^1(F^*))$ is the inclusion $V^*\hookrightarrow H^1(F^*)$. Moreover, we have:
$$\begin{array}{rcl} rank(E)& = & rank(F)+ \dim V= r+v  , \\
c_1(E) & = & c_1(F) \\
c_2(E)& = & c_2(F).
\end{array}
$$
\end{definition}
We define $U_e\subset  \PSpl(r; L,c_2)$ as the open locus parameterizing  vector bundles $F$ on $X$ such that $H^0(F)=H^0(F^*)=0$ and we assume $U_e\ne \emptyset $. We fix an integer $v\ge 1$ and we consider the 
Grassmann bundle $\pi: \cP_{U_e}\to U_e$ where 
$$\cP _{U_e}:=\{(F,V) \mid F\in U_e, V^*\subset H^1(F^*), \dim (V)=v\}.$$
$\cP_{U_e}$ is a smooth algebraic space of dimension $$\dim \cP_{U_e}=\dim \PSpl(r;L,c_2)+ v(u-v)$$ where $u:=\dim H^1(F^*)=-\chi (F^*)$.
By \cite[Proposition 3.6]{FM2} there is a natural morphism

$$  \beta: \cP_{U_e}   \to   \PSpl(r+v;L,c_2) 
$$ 
 which  extends the set theoretic map $(F ,V)$$  \mapsto  E$ where $E$ is the extension bundle associated to $(F,V)$.
The morphism $\beta $ is an injective locally closed embedding and its differential map is also injective. 

\begin{remark}
    In \cite{FM2} we do not fix the determinant; however it is immediate by the definition that there is a canonical isomorphism $\det F\to \det E$ which extends to families.
\end{remark}

For any $\cL\subset \PSpl(r;L,c_2)$ smooth irreducible locally closed subspace, we define $\cP _{\cL}:=\pi ^{-1}(\cL\cap U_e)$ and we will denote by $\widetilde{\cL}\subset \PSpl(r+v;L,c_2)$ the image of $  \beta_{|\cP_{\widetilde{\cL}}}: \cP_{\widetilde{\cL}}   \to   \PSpl(r+v;L,c_2).$

\vskip 6mm
\noindent {\bf Proof of Theorem 1.2.} It immediately follows from \cite[Theorem 5.7]{FM2} since the hypothesis that $X$ is smooth was never used. \qed

\subsection{Generalized syzygy bundles} We start the subsection quickly recalling the construction of generalized syzygy bundles and for more details the reader can look at \cite{FM1} and \cite{FM2}.

\begin{definition}\label{def1}
Let $F$ be a globally generated vector bundle of rank $r$, determinant $L$ and second Chern class $c_2$ on a K3  surface $X$. For any general subspace $W\subset  H^0(X,F)$ with $w:=\dim W\ge 2+r$ we define the {\em generalized syzygy bundle} $S$ on $X$ as the kernel of the evaluation map $$eval_{W}: W\otimes \mathcal O_X \to F.$$ Therefore, $S$  fits inside a short exact sequence:
\begin{equation}
    \label{exact1}
e: \quad 0\to  S\to W\otimes \mathcal O_X \to  F\to 0, \text{ and }
\end{equation}
$$\begin{array}{rcl} rank(S)& = & \dim W- rank(F)= w-r, \\
c_1(S) & = & -c_1(F) \\
c_2(S)& = & -c_2(F)+c_1(F)^2.
\end{array}
$$
\end{definition}

We will denote by $$U:=U(r;L,c_2)\subset  \PSpl(r; L,c_2)$$ the open locus parameterizing  globally generated rank $r$ perfect simple vector bundles $F$ on $X$ such that $H^1(\cF)=0$ and we define
$$\cG _U:=\{(F,W) \mid F\in U, W\subset H^0(F), \dim (W)=w\}.$$
The natural projection $\pi: \cG_U\to U$ is a Grassmann bundle and  $\cG_U$ is a smooth algebraic space of dimension $$\dim \cG_U=\dim \PSpl(r;L,c_2)+ w(v-w)$$ where $v:=\dim H^0(F)$.
Moreover, we define 
$$\cG^0 _U:=\{(F,W)\in \cG _U  \mid eval_W \text{ is surjective}\}$$
which is open in $\cG _U$ and we keep denoting  by $\pi $ the natural projection $\pi : \cG^0_U \to U$.  For more details, see \cite[Definition 3.4]{FM1} and \cite[Definition 3.7]{FM2}.
 Arguing as in \cite[Proposition 3.5]{FM1} we check that for any K3 surface $X$, the generalized syzygy bundle associated to any pair $(F,W)\in \cG^0_U$ is simple and perfect. Moreover, there is a natural morphism

$$  \alpha: \cG_U^0   \to   \PSpl(w-r;-L,L^2-c_2) 
$$ 
 which  extends the set theoretic map $$(F ,W)  \mapsto  S:= ker(  W\otimes \cO_X \twoheadrightarrow F).$$
The morphism $\alpha $ is an injective locally closed embedding and its differential map is also injective.

For any $\cL\subset \PSpl(r;L,c_2)$ smooth irreducible locally closed subspace, we define $\cG^0 _{\cL}:=\pi ^{-1}(\cL\cap U)$ and we will denote by $\cL'\subset \PSpl(w-r;-L,L^2-c_2)$ the image of $  \alpha_{|\cG_{\cL}^0}: \cG_{\cL}^0   \to   \PSpl(w-r;-L,L^2-c_2).$
\vskip 6mm
\noindent {\bf Proof of Theorem 1.3.} (1) The fact that $$\dim \cL = \frac{1}{2}\dim Ext^1(F,F)_0=\frac{1}{2}\dim \PSpl(r;L,c_2)$$ if and only if $$\dim \cL' = \frac{1}{2}\dim Ext^1(S,S)_0=\frac{1}{2}\dim \PSpl(w-r;-L,L^2-c_2)$$ follows from \cite[Proposition 3.19]{FM2} where the hypothesis of K3 being smooth is not necessary.

(2) We consider a smoothing $(\pi :\cX\longrightarrow B,\cL)$ of the pair $(X,L)$, the  relative moduli spaces $f: \cM=\PSpl_{\pi}(r;\cL,c_2)\to B$ and $f': \cM'=\PSpl_{\pi}(w-r;-\cL,\cL^2-c_2)\to B$ and the open locus $Q\subset Quot(\cO_{\cX}^{\oplus w}/B)$ parametrizing $[0\to S \to \cO_{\cX_b}^{\oplus w}\to F\to 0]$ such that $F$ and $S$ are simple vector bundles on $\cX_b$ of rank and Chern classes $(r;\cL_{|\cX_b},c_2)$ and $(w-r,-\cL_{|\cX_b}, \cL_{|\cX_b}^2-c_2)$, respectively. 

We have the commutative diagram:
\[
\xymatrix{ Q \ar[r]^{g} \ar[d]^h & \cM \ar[d]^{f }\\
\cM'\ar[r]^{f'} & B.
}
\]
For any $\alpha \in H^0(\cM,\Omega ^2_{\cM/B})$ and $\beta  \in H^0(\cM',\Omega ^2_{\cM'/B})$ the 2-form $g^*\alpha +h^*\beta\in H^0(Q,\Omega ^2_{Q/B})$ is zero if and only if it is zero pointwise. By \cite[Proposition 4.3]{FM2}, we know that it is zero  for any $b\in B\setminus \{b_0\}$; therefore it is zero everywhere. 
    \qed

\section{Smoothability criteria for polarized K3 surfaces}\label{criteriasmoothability}

\subsection{Smoothability for polarized schemes}

Throughout this subsection, we assume that $X$ is a reduced, connected, pure dimensional, projective scheme over $\CC$ with lci singularities (i.e., for some/any closed embedding $i:X\to M$ in a smooth projective variety $M$, the closed embedding $i$ is regular).

This implies that $\EExt^i(\Omega_X,\cO_X)=0$ for all $i>1$. We denote by $T_X$ the tangent sheaf $\Omega_X^\vee$ and by $T^1_X$ the sheaf $\EExt^1(\Omega_X,\cO_X)$. More generally, for any flat morphism $f:\cX\to B$ with lci reduced fibers we have $\EExt^i(\Omega_f,\cO_\cX)=0$ for all $i>1$, and we denote by $T_f$ the tangent sheaf $\Omega_f^\vee$ and by $T^1_f$ the sheaf $\EExt^1(\Omega_f,\cO_\cX)$. The assumption that fibers are lci implies that $T_f$ and $T^1_f$ commute with base change \cite[Theorem 2.11]{FFP2}.

\begin{definition}
    A {\em formal one-parameter deformation} of $X$ over the formal spectrum of $\CC[[t]]$ is the datum of \begin{enumerate}
    \item for every $n\ge 0$, a flat morphism $\pi_n:X_n\to S_n:=\Spec \CC[t]/t^{n+1}$;
    \item an isomorphism $\phi:X\to X_0$;
    \item for  every $n\ge 0$, a morphism $f_n:X_n\to X_{n+1}$ 
\end{enumerate}  
such that for every $n\ge 0$ the induced  diagram 
\[
\xymatrix{ X_n \ar[r]^{f_n} \ar[d]^{\pi_n} & X_{n+1} \ar[d]^{\pi_{n+1}}\\
S_n\ar[r]^{i_n} & S_{n+1}
}
\]
 is commutative and cartesian, where $i_n:S_n\to S_{n+1}$ is the closed embedding such that $i_n^*t=t$. Note that this makes each $\pi_n$ into an infinitesimal deformation of $X$ over $S_n$.
\end{definition}
In \cite[Definition 11.6]{Tz}, Tziolas defines the notion of formal smoothing as a formal one-parameter deformation satisfying a number of equivalent conditions. 
In this paper, we will use as definition the characterization of formal smoothing for lci reduced schemes given by Nobile in \cite[Theorem 5.7]{N}.
\begin{definition}\label{Def-Tz}
    A formal one-parameter deformation of $X$ is called a {\em formal smoothing} if there exists an $N>0$ such that the natural map $f_{N*}T^1_{\pi_N}\to T^1_{\pi_{N+1}}$ (induced by the isomorphism $f_N^*T^1_{\pi_{N+1}}\cong T^1_{\pi_N}$) is an isomorphism. 
    
    Note that by Nakayama's lemma, we have isomorphisms $f_{n*}T^1_{\pi_n}\to T^1_{\pi_{n+1}}$ for every $n\ge N$.
\end{definition}

The key property of formal smoothings is the following Proposition, due to Tziolas (\cite[Proposition 11.8]{Tz}. 
\begin{proposition}\label{smoothing_is_formal}
    Assume $f:\cX\to B$ is a flat projective morphism with base a smooth curve $B$. Let $b_0\in B$, choose $t\in \cO_{B,b_0}$ a local parameter, and let $X:=\cX_{b_0}$. Then $f$ induces a formal one-parameter deformation of $X$.
If $X$ is reduced and lci, then the general fiber of $f$ is smooth (i.e., $f$ is a smoothing) if and only if the associated formal deformation is a formal smoothing.
\end{proposition} 

Tziolas gave the following sufficient conditions for the existence of a formal smoothing of $X$: 
\begin{theorem} 
    Assume that $H^2(X,T_X)=H^1(X,T^1_X)=0$, and that $T^1_X$ is generated by global sections. Then $X$ admits a formal smoothing.
\end{theorem}
\begin{proof}
See \cite[Theorem 12.5]{Tz}.
\end{proof}

\begin{definition}
    Let  $(X,L)$ be a polarized scheme: we assume that $X$ is reduced, projective and lci, and that $L$ is ample on $X$.
A {\em formal one-parameter polarized deformation} of $(X,L)$ is the data 
of \begin{enumerate}
    \item a formal one-parameter deformation $(X_n,f_n,\pi_n,\phi)$ of $X$;
    \item for every $n\ge 0$, a line bundle $L_n\in \Pic(X_n)$ and an isomorphism $\alpha_n:f_n^*L_{n+1}\to L_n$;
    \item an isomorphism $\alpha:\phi^*L_0\to L$.
\end{enumerate}
We say that $(X,L)$ is {\em formal polarised smoothable} if there exists a formal one-parameter polarized deformation of $(X,L)$ such that the underlying formal one-parameter deformation of $X$ is a formal smoothing.
\end{definition}

\begin{proposition}\label{formal_pol_smoothable}
If $(X,L)$ is a polarized scheme (we always assume $X$ reduced with lci singularities) which is formal polarized smoothable  then it is  smoothable in the sense of Definition \ref{def_smoothable}.    
\end{proposition}
\begin{proof}
   By replacing $L$ with $L^{\otimes N}$ with $N\gg 0$, we may assume that $L$ is very ample on $X$ and that $H^1(X,L)=0$. Let $m+1:=h^0(X,L)$. From this it follows that $L_n$ is very ample for $X_n$ relative to $S_n$, and that $R^1\pi_{n*}(L_n)=0$ by cohomology and base change. By the same theorem, $E_n:=\pi_{n*}(L_n)$ commutes with base change, and it is locally free of rank $m+1$ since $H^{-1}(X,L)=0$; $E_n$ must be globally trivial since $S_n$ has only one point.
   
By  definition \ref{Def-Tz}, we can choose an integer $N$ such that the fact that $\{X_n\}$ is a formal smoothing can be checked at level $N$. We choose a trivialization of $E_{N}$, which induces a closed embedding $X_{N}\to \mathbb P^m_{S_{N}}$, hence a morphism $\varphi :S_{N}\to \Hilb^P(\mathbb P^m)$ where $P$ is the Hilbert polynomial of $(X,L)$.
 We have constructed a map $S_{N}\to \Hilb^P(\mathbb P^m)$ which can be extended to $S_n$ for every $n\ge N$ by lifting the trivialization of $E_{N}$ to a trivialization of $E_n$.

 By Artin's approximation theorem \cite{A1}, we can find \begin{enumerate} 
 \item a morphism $\psi:B\to \Hilb^P(\mathbb P^m)$ where $B$ is a smooth affine curve,
 \item a point $b\in B$,
 \item a generator $u$ of $m_b\subset \cO_{B,b}$,
 \end{enumerate}
 such that $\psi\circ \rho=\varphi$ where, $\rho:S_{N}\to B$ is the unique morphism mapping $S_{N,red}$ to $b$ and such that  $\rho^*u=t$.
 Let $X_B\to B$ be the flat family corresponding to $\psi$, and consider now the formal deformation $(X_n',L_n')$ of $X$ defined by $X_n':=X_B\times_BS_n$ and $L_n'$ the pullback of $\cO_{X_B}(1)$
 for every $n$. This is a formal deformation of $(X,L)$ with  $(X_n',L_n')=(X_n,L_n)$
  for every $n\le N$. By definition of $N$, 
  it must be a formal smoothing. By Proposition \ref{smoothing_is_formal}, the flat family $X_B\to B$ is a smoothing, hence  $(X_B\to B, \cO_{X_B}(1))$ is a smoothing of the pair $(X,L)$ in the sense of Definition \ref{def_smoothable}.
  \end{proof}

\begin{remark}
    A similar argument in a simpler context, where the ample line bundle is the dualizing sheaf and hence canonically extends to every infinitesimal deformation, can be found in \cite{FFP}.
\end{remark}

\subsection{Smoothability of polarized K3 surfaces with isolated singularities}

In what follows,  $X$ is a singular reduced $K3$ surface with isolated lci singularities, and $L\in \Pic(X)$ an ample line bundle. 
The aim of this section is to generalize  Tziolas' criterion for formal  smoothability to the pair $(X,L)$. Since $X$ has isolated singularities, Tziolas' criterion reduces to Burns and Wahl \cite{BW}, namely $H^2(X,T_X)=0$.

We recall that we can associate to $L$ its Atiyah sequence \[
0\to \Omega_X\to \cP\to \cO_X\to 0
\]
which defines a class $e(L)\in H^1(\Omega_X)$; since $\cO_X$ is locally free, the dual sequence is also exact 
\[
0\to \cO_X\to \cP^\vee \to T_X\to 0
\]
and we have an induced isomorphisms $\EExt^1(\cP,\cO_X)\cong \EExt^1(\Omega_X,\cO _X)= T^1_X$ while $\EExt^i(\cP,\cO_X)=0$ for every $i>1$.

\begin{lemma}
    Assume that $H^2(X,T_X)=0$; then $\Ext^2(\Omega_X,\cO_X)$ is also zero.
\end{lemma} 
\begin{proof}
    It is enough to show that the other two terms of the local to global spectral sequence also vanish: $H^1(X,T^1_X)$ because $X$ has isolated singularities, and $H^0(X,\EExt^2(\Omega_X,\cO_X))$ because since $X$ is lci then the sheaf $\EExt^2(\Omega_X,\cO_X)$ is zero.
\end{proof}
In view of this Lemma, under the assumption that $H^2(X,T_X)=0$ the short exact sequences above and the local-to-global spectral sequence of Ext yield a commutative diagram with exact rows and columns

\[
\begin{tikzcd}
0\arrow{d} & 0\arrow{d}\\
     H^1(\cP^\vee) \arrow[hookrightarrow]{r}\arrow{d} & \Ext^1(\cP,\cO_X) \arrow{r}\arrow{d} & H^0(T_X^1) \arrow[equal]{d}\arrow[twoheadrightarrow]{r} & W\arrow{d}\\
  H^1(T_X) \arrow[hookrightarrow]{r}\arrow{d} & \Ext^1(\Omega_X,\cO_X) \arrow{r}\arrow{d} & H^0(T_X^1) \arrow{r} & 0\\
    H^2(\cO_X)\arrow[equal]{r} \arrow{d} &H^2(\cO_X)\arrow{d}\\
    H^2(\cP^\vee)\arrow{r}\arrow{d} &\Ext^2(\cP,\cO_X)\arrow{d}\arrow{r} & 0\\
    0 &  0
\end{tikzcd}
\]
where $W$ fits into the exact sequence \[
0\to W\to H^2(\cP^\vee) \to \Ext^2(\cP,\cO_X) \to 0.
\]
\begin{remark}
    We briefly recall the geometric meaning of various groups and maps in the diagram above. $H^1(\cP^\vee)$ and $H^2(\cP^\vee)$ are tangent and obstruction spaces to locally trivial deformations of the pair $(X,L)$; $\Ext^1(\cP,\cO_X)$ and $\Ext^2(\cP,\cO_X)$ are tangent and obstruction spaces to all deformations of the pair; $H^1(T_X)$ and $H^2(T_X)$ (resp.~$\Ext^1(\Omega_X,\cO_X)$ and $\Ext^2(\Omega_X,\cO_X)$) are tangent and obstruction spaces to locally trivial (respectively all) deformations of $X$; $H^1(\cO_X)$ and $H^2(\cO_X)$ are relative tangent and obstruction to deformations of the pair $(X,L)$ over deformations of $X$, and thus also for the locally trivial case. Note that so far we never used that the singularities are isolated other than to guarantee that $H^1(X,T^1_X)=0$.
    
    Since $X$ has isolated singularities, the vector space $H^0(T^1_X)$ is the product, over the singular points $p$ in $X$, of the tangent spaces to the deformation functor $D_{X,p}$ of the germ of $p$ at $X$; since $X$ is lci, all local deformation functors $D_{X,p}$ are unobstructed. 
\end{remark}
Recall that since $H^2(\cO_X)$ is one dimensional, any complex linear map to it is either zero or surjective.
\begin{lemma}\label{surjective}
     The morphism $\Ext^1(\Omega_X,\cO_X)\to H^2(\cO_X)$ in the diagram above is surjective for any singular K3 surface $X$.
\end{lemma}
\begin{proof}
   This morphism is induced by the cup product \[
    \Ext^1(\cO_X,\Omega_X)\otimes\Ext^1(\Omega_X,\cO_X)\to \Ext^2(\cO_X,\cO_X)=H^2(\cO_X)
    \]
    applied to the Atiyah class of $L$ in $\Ext^1(\cO_X,\Omega_X)=\Ext^1(L,\Omega_X\otimes L)$, which corresponds to the Atiyah sequence for $L$\[
    0\to \Omega_X\otimes L\to \cP_L\to L\to 0.
    \]
    It is not zero because the Atiyah class of $L$ is nonzero since $L$ is ample.
\end{proof}

\begin{corollary}\label{obstructions}
    For any singular  K3 surface $X$, the natural morphism $\Ext^2(\cP,\cO_X)\to \Ext^2(\Omega_X,\cO_X)$ is an isomorphism.
    \begin{proof}
        From the exact sequence defining $\cP$ we have that the cokernel is contained in $H^3(\cO_X)=0$, and the kernel is equal to the cokernel of
        $\Ext^1(\Omega_X,\cO_X)\to H^2(\cO_X)$, which is zero by the Lemma \ref{surjective}.
    \end{proof}
\end{corollary}
\begin{theorem}\label{main_smoothing} Let $X$ be a singular K3 surface with isolated lci singularities, and $L\in \Pic(X)$ an ample line bundle. The pair $(X,L)$ is polarized smoothable
    if \begin{enumerate} \item $H^2(X,T_X)=0$, and 
    \item the morphism  $\psi:\Ext^1(\cP,\cO_X) \to H^0(T_X^1)$ is surjective.
\end{enumerate}

\end{theorem} 

\begin{proof}
    From Lemma \ref{surjective} it follows that  $\Ext^2(\cP,\cO_X)=0$. 

    It is therefore enough to show that if $\Ext^2(\cP,\cO_X)=0$ and $\psi$ is surjective, then $(X,L)$ is polarized smoothable. Since $\Ext^2(\cP,\cO_X)=0$, the deformation functor of the pair $(X,L)$ is unobstructed. Hence any infinitesimal deformation over any $S_n$ can be extended to a formal one and it is enough to find a deformation of $(X,L)$ over some $S_{N+1}$ that satisfies the definition of formal smoothing (i.e., $f_{N*}T^1_{\pi_N}\to T^1_{\pi_{N+1}}$ is an isomorphism), and apply Proposition \ref{formal_pol_smoothable}. 

    Let $p_1,\ldots,p_m$ be the singular points of $X$. To each of them we  associate the infinitesimal local deformation functor $D_{(X,p_i)}$,
    which is unobstructed (because $X$ is lci) and has as tangent space the stalk $T^1_{X,p_i}$ of $T^1_X$ at $p_i$, which is a finite dimensional vector space; the natural map \[
    H^0(X,T_X^1)\to \bigoplus_{i=1}^m T^1_{X,p_i}
    \]
    is an isomorphism.
    
    Each isolated lci singularity $(X,p_i)$ is smoothable. 
     So we can find a local formal deformation $(\alpha_{i,n})_{n\in \NN}$ which
     is a formal smoothing. Let $D_{loc}:=\prod D_{X,p_i}$ 
     be the local deformation functor of $X$, and for $n\ge 0$ let $\alpha_n\in D(S_n)$ be the product of all  $\alpha_{i,n}$. 

   By the assumption (2), the morphism of functors $D_{(X,L)}\to D_{loc}$ induces a surjection on tangent spaces and it is smooth since both functors are unobstructed, as assumption (1) implies that the obstruction space $\Ext^2(\cP,\cO_X)$ of $D_{(X,L)}$ is zero. Hence the sequence $\alpha_n$ can be lifted to a formal deformation $(X_n,L_n)$ of the pair $(X,L)$, which is a formal smoothing since each $\alpha_{i,n}$ smooths the point $p_i$. By definition, there exists an $N>0$ such that $f_{N*}T^1_{\pi_N}\to T^1_{\pi_{N+1}}$ is an isomorphism.

    We conclude the proof by applying Proposition \ref{formal_pol_smoothable}.
\end{proof}
\begin{remark}
Assumption (2) of the Theorem \ref{main_smoothing} is verified whenever $H^1(X,T_X)\to H^2(X,\cO_X)$ is surjective. In this case both morphisms $H^1(X,\cP^\vee)\to H^1(X,T_X)$ and $\Ext^1(\cP,\cO_X) \to \Ext^1(\Omega_X,\cO_X)$ are codimension one embeddings. Hence by counting dimension $W=0$ and $\Ext^1(\cP,\cO_X)\to H^0(T^1_X)$ is surjective.
\end{remark}

\subsection{The case of quartics in $\PP^3$}

We end this paper with a concrete example which illustrates very well how our criterion (Theorem \ref{main_smoothing}) applies, and moreover that none of the hypotheses  of Theorem \ref{main_smoothing} is necessary. From now on, we will assume that $X\subset \PP^3$ is a quartic surface with isolated singularities (necessarily hypersurface singularities, in particular lci). 

The pair $(X,\cO_X(1))$ is clearly polarized smoothable; the deformations of the pair are unobstructed and their tangent space $\Ext^1(\cP,\cO_X)$ has dimension $h^0(\PP^3, \cO_{\PP^3}(4))-\dim GL(4)=19$ because all such deformations are quartics in $\PP^3$.

\begin{lemma}\label{key24}
In the assumption that $X$ is a quartic  surface in $\PP^3$ with isolated singularities, we have\begin{enumerate}
    \item $\Ext^2(\cP,\cO_X)=\Ext^2(\Omega_X,\cO_X)=0$;
    \item $H^0(X,T_X)=0$;
    \item $H^2(X,T_X)=0$ if and only if $\Ext^1(\Omega_X,\cO_X)\to H^0(T^1_X)$ is surjective;
    \item  $\Ext^1(\cP,\cO_X)\to H^0(T^1_X)$ is surjective if and only if $H^1(\PP^3,\cJ(4))=0$ where $\cJ $ is the Jacobian ideal sheaf of $X$.
\end{enumerate}    
\end{lemma}
\begin{proof}
(1) From the short exact sequence \[
0\to \cO_X(-4)\to \Omega_{\PP^3}|_X\to \Omega_X\to 0
\] 
we get the exact sequence 
\[
H^1(X,\cO_X(4))\to \Ext^2(\Omega_X,\cO_X)\to H^2(T_{\PP^3}|_X).
\]
We have $H^1(X,\cO_X(4))=0$ by Kodaira vanishing. By the Euler exact sequence $H^2(T_{\PP^3}|_X)$ is a quotient of $H^2(\cO_X(1))^{\oplus 4}$ and hence is also zero, again by Kodaira. Thus $\Ext^2(\Omega_X,\cO_X)=0$, which implies $\Ext^2(\cP,\cO_X)=0$ by Corollary \ref{obstructions}. \\
(2) $H^0(T_X)=0$ because it is Serre dual to $\Ext^2(\Omega_X,\cO_X)$. \\
(3) From the long exact sequence of $\Ext$, we get that $H^2(X,T_X)=0$ if and only if $\Ext^1(\Omega_X,\cO_X)\to H^0(T^1_X)$ is surjective. In particular this is true if $\Ext^1(\cP,\cO_X)\to H^0(T^1_X)$ is surjective. \\
(4) Under our assumptions, $T^1_X=\cO_Z(4)$ where $Z$ is the scheme theoretic support of $T^1_X$, i.e. the scheme structure on $\Sing(X)$ induced by the Jacobian ideal sheaf $\cJ\subset \cO_{\PP^3}$.

Thus, by the exact sequence above we have an exact sequence \[
0  \to H^0(T_{\PP^3}|_X)\to H^0(\cO_X(4))\to \Ext^1(\Omega_X,\cO_X)\to H^1(T_{\PP^3}|_X)\to 0
\]
and by the Euler sequence $H^1(T_{\PP^3}|_X)=H^2(X,\cO_X)\cong \CC$. Hence we get an exact sequence 
\[
0 \to H^0(T_{\PP^3}|_X)\to H^0(\cO_X(4))\to \Ext^1(\cP,\cO_X)\to 0.
\]

Therefore, the map $\Ext^1(\cP,\cO_X)\to H^0(T^1_X)$ is surjective if and only if the restriction map $H^0(\cO_X(4))\to H^0(\cO_Z(4))$ is surjective. Since $H^1(\cO_X(4))=0$ by Kodaira, the cokernel of this map is isomorphic to $H^1(\cJ(4))$.
\end{proof}

Whether $H^2(X,T_X)$ or $H^1(\PP^3,\cJ(4))$ are zero can be checked very explicitly if one knows the singularities of $X$.

\vskip 2mm
We now give three explicit examples where we can see which of the two assumptions of the Theorem \ref{main_smoothing} are satisfied.  In Example \ref{ex1} we have that assumption (2) of Theorem \ref{main_smoothing} fails (we do not know if (1) also fails);  in Example \ref{ex2} both assumptions (1) and (2) of Theorem 5.12 fail.
Finally, in Example \ref{ex3} both assumptions (1) and (2) of Theorem 5.12 are satisfied. Note that if assumption (2) holds, then so does (1).

\vskip 4mm

\begin{example}  \label{ex1}
     Let $X\subset \PP^3=\Proj(\CC[x,y,z,t])$ be the quartic surface defined by  $f=xy^3+yz^3+t^4$ and let   $J=(t^3, y^3, 3xy^2+z^3, yz^2)\subset \CC[x,y,z,t]$ be its jacobian ideal. We easily check that  $\Sing (X)=\{(1,0,0,0)\}$.  Therefore, $T^1_X=\cO_Z(4)$ where $Z$ is the scheme  structure on $\Sing(X)$ induced by the Jacobian ideal sheaf $\cJ=\widetilde{J}\subset \cO_{\PP^3}$. We notice that $J$ is a saturated ideal with the following locally free resolution:\[
   0\to \cO(-7)\oplus \cO(-8)\xrightarrow{\,\,\,C\,\,\,} 
   \cO(-4)\oplus \cO(-5)\oplus \cO(-6)^{\oplus 3}\xrightarrow{\,\,\,B\,\,} 
   \cO(-3)^{\oplus 4}\xrightarrow{\,\,\,A\,\,\,}  \cJ\to 0
    \]
  where 
  \[ C:=\begin{pmatrix} t^3 & 0 \\
  0 & t^3 \\
  -y & 0 \\
  z/3 & y^2      
  \end{pmatrix}, \quad
  B:=\begin{pmatrix} -y & 0 & -t^3 & 0 & 0 \\
  x & -z^2 & 0 & -t^3 & 0 \\
  z/3 & y^2 & 0 & 0 & -t^3 \\
  0 & 0 & xy^2+z^3/3 & y^3 & yz^2 
  \end{pmatrix}, \quad \text{and}
  \]
  \[
  A:=\begin{pmatrix} 
      xy^2+z^3/3 & y^3 & yz^2 & t ^3
  \end{pmatrix}.
   \]
    We split the last exact sequence into two short exact sequences and twist them by $\cO(4)$:
    \[
    0\to \cO(-3)\oplus \cO(-4)\to
   \cO\oplus \cO(-1)\oplus \cO(-2)^{\oplus 3}\to \cF(4)\to 0
    \]
and\[
0\to\cF(4)\to \cO(1)^{\oplus 4}\to \cJ(4)\to 0.
\]
Taking cohomology  we get \[
 H^1(\cJ(4))\cong H^2(\cF(4))\cong  H^3(\cO(-3)\oplus \cO(-4))\cong \CC.
\]
From Lemma \ref{key24} (4) we conclude that $\Ext^1(\cP,\cO_X)\to H^0(T^1_X)$ is  not surjective.
\end{example}

\begin{example}  \label{ex2}
  Let $X\subset \PP^3=\Proj(\CC[x,y,z,t])$ be the quartic surface defined by  $f=t^4+xy(x-y)(x+y)$ and let   $J=(3x^2y-y^3,x^3-3xy^2,t^3)\subset \CC[x,y,z,t]$ be its jacobian ideal. We easily check that  $\Sing (X)=\{(0,0,1,0)\}$.  Therefore, $T^1_X=\cO_Z(4)$ where $Z$ is the scheme  structure on $\Sing(X)$ induced by the Jacobian ideal sheaf $\cJ=\widetilde{J}\subset \cO_{\PP^3}$. Notice that $J$ is a saturated ideal, it is a complete intersection ideal and it
  has the locally free resolution given by the Koszul complex:
\[
   0\to \cO(-9)\to 
   \cO(-6)^{\oplus 3}\to 
   \cO(-3)^{\oplus 3}\to  \cJ\to 0.
    \]
 Arguing as in the previous example shows that $H^1(\cJ (4))\cong H^3(\cO(-5))$ which has dimension $4$, thus $\Ext^1(\cP,\cO_X)\to H^0(T^1_X)$ is  not surjective  (Lemma \ref{key24} (4)). Finally, we consider  the exact sequence
 \[
0\to T_X \to T_{\PP^3}|_X \to \cO _X(4) \to T^1_X\cong \cO _Z (4)\to 0
 \]
we cut it into two short exact sequences, we take cohomology and we get $\dim H^2(X,T_X)\ge 3$ using that  $H^2(\PP^3,\cJ (4))\cong H^2(X,\cJ_{Z/X}(4))\cong \CC^4$ and $H^1(T_{\PP^3}|_X)\cong \CC$. 

This example extends, by the same argument, to any projective cone over a smooth quartic curve.
\end{example}

\begin{example}  \label{ex3}
  Let $X\subset \PP^3=\Proj(\CC[x,y,z,t])$ be the quartic surface defined by  $f=t^4+x^2y^2+x^2z^2+y^2z^2$, Let   $$J=(x(y^2+z^2),y(x^2+z^2),z(x^2+y^2),t^3 )\subset \CC[x,y,z,t]$$ be its Jacobian ideal and $J^{sat}=(xy,xz,yz,t^3)$ the saturated ideal. 
  \\
  We have $\Sing (X)=\{(1,0,0,0),(0,1,0,0),(0,0,1,0)\}$ and $T^1_X=\cO_Z(4)$ where $Z$ is the scheme  structure on $\Sing(X)$ induced by the Jacobian ideal sheaf $\cJ=\widetilde{J^{sat}}\subset \cO_{\PP^3}$. We easily check that $\cJ$ 
  has the followng locally free resolution:
\[
   0\to \cO(-6)^{\oplus 2}\to 
   \cO(-5)^{\oplus 3}\oplus \cO(-3)^{\oplus 2}\to 
   \cO(-2)^{\oplus 3}\oplus \cO(-3)\to  \cJ\to 0.
    \]
 We immediately get that $H^1(\PP^3,\cJ (4))=0$ and hence  $\Ext^1(\cP,\cO_X)\to H^0(T^1_X)$ is   surjective, thus $\Ext^1(\Omega_X,\cO_X)\to H^0(T^1_X)$  and $H^2(X,T_X)=0$.
\end{example}

We remark that the morphism $X\to \PP^2$ given by $(x,y,z,t)\mapsto (x,y,z)$ makes $X$ into a simple cyclic cover of order $4$ of the plane branched over the quartic curve $x^2y^2+x^2z^2+y^2z^2$, which has three nodes and no other singularity. Hence $X$ has three $A_3$ singularities. The explicit calculation
above is a special case of the following more general result
\begin{proposition}
    Let $(X,L)$ be a singular projective K3 surface with only ADE singularities; then $\Ext^1(\cP,\cO_X)\to H^0(T^1_X)$ is   surjective, hence $\Ext^1(\Omega_X,\cO_X)\to H^0(T^1_X)$ is also surjective and $H^2(X,T_X)=0$. 
\end{proposition}
\begin{proof}
    Let $\eps:\tilde X \to X$ be a minimal resolution of singularities; then $\tilde X$ is a nonsingular projective $K3$ surface. Let $\{E_i\}_{i\in I}$ be the irreducible components of the exceptional divisors of $\eps$; they are smooth rational curves with $E_i^2=-2$. Let $e_i$ be the class of $E_i$ in $H^1(\tilde X,\Omega^1_{\tilde X})\subset H^2(\tilde {X},\QQ)$. Over each singular point the intersection form $(e_ie_j)$ is negative definite, thus the $e_i$ are linearly independent. Let $h:=c_1(\eps^*L)$; since $h\cdot e_i=0$ for each $i\in I$, it follows that also $(h,e_i)_{i\in I}$ are linearly independent. This implies that the natural map $H^1(T_{\tilde X})\to H^2(\cO_{\tilde X})\oplus H^2(\cO_{\tilde X})^{I}$ given by cup product with $(h,e_i)_{i\in I}$ is surjective.

Let $U\subset X$ be the smooth open locus, and for each singular point $p$, choose an affine open $V_p$ such that $V_p\cap \Sing(X)=\{p\}$. Let $\tilde U:=\eps^{-1}(U)$ and $\tilde V_p:=\eps^{-1}(V_p)$. Let $I_p:=\{i\in I\,|\, E_i\subset \tilde V_p\}$.

The morphism $\eps:\tilde U\to U$ is an isomorphism, inducing an isomorphism $D_{\tilde U}\to D_U$. By \cite{Br} there is a natural action of the Weyl group $G_p$ of the ADE singularity at $p$ on $D_{\tilde V_p}$ such that the quotient is naturally isomorphic to $D_{V_p}$. Both $D_{\tilde V_p}$ and $D_{V_p}$ are unobstructed; the action of $G_p$ on $T^1D_{\tilde V_p}=\bigoplus_{i\in I_p}H^1(E_i,N_{E_i})$ is induced by a permutation of the $E_i$. Note that the natural map $H^1(E_i,N_{E_i})\to H^2(\cO_{\tilde X})$ is an isomorphism, and it maps $H^1(\tilde X, T_{\tilde X})\to H^1(E_i,N_{E_i})$ to the cup product with the class of $E_i$ in $H^1(\Omega^1_{\tilde X})$.

It follows that $D_{(\tilde X, \eps^*L)}\to \prod_pD_{\tilde V_p}$ is a smooth functor. The action of $\prod_p G_p$ on $\prod_p D_{\tilde V_p}$ lifts to $D_{(\tilde X,\eps^*L)}$ by viewing it as the fiber product of all $D_{\tilde V_p}$ and $D_{\tilde U}$ and $D_{(X, L)}$ over the product of all $D_{V_p}$ and $D_U$.

Taking the quotient by $\prod_pG_p$ shows that $D_{(X,L)}\to \prod_pD_{V_p}$ is also surjective on tangent spaces, which is the claim.
\end{proof}


\begin{thebibliography}{99}


\bibitem{Ar} I.V.~Artamkin,  {\em On deformation of sheaves}, Mathematics of the USSR-Izvestiya {\bf 32(3)} (1989), 663--668.

\bibitem{A1} M. F.~Artin, {\em Algebraic approximation of structures over complete local rings},
Publications math\'ematiques de l'IHES {\bf 36} (1969), 23-58.

\bibitem{A2}  M. F.~Artin, {\em Algebraic construction of Brieskorn's resolutions}, Journal of Algebra {\bf 29} (1974),  330--348.

\bibitem{BHPV} W.P.~Barth, K.~Hulek, C.A.M.~Peters, A.~van de Ven, {\em Compact Complex Surfaces}, Springer Verlag, 2004.

\bibitem{Br} E.~Brieskorn, {\em Singular elements of semisimple algebraic groups}, in {\em Actes Congr\`es Intern.\ Math.\ 1970}, 279--284.

\bibitem{BW}
D.M.~Burns Jr. and J.M.~Wahl {\em  Local contributions to global deformations of surfaces},
Inventiones mathematicae {\bf 26} (1974), 67-?88.

\bibitem{FM1} B.~Fantechi and R. M.~Mir\'o-Roig, {\em  Moduli of generalized syzygy bundles}, Preprint, arXiv 2306.04317v1.  Pure and Applied Mathematics Quarterly, to appear.


\bibitem{FM2} B.~Fantechi and R. M.~Mir\'o-Roig, {\em Lagrangian subspaces of the moduli space of simple sheaves on K3 surfaces}, Preprint arXiv:2306.05338.

\bibitem{FFP} B.~Fantechi, M.~Franciosi, R.~Pardini, {\em Smoothing semi-smooth stable Godeaux surfaces},  Algebraic Geometry {\bf 9} (2022) 502-512.

\bibitem{FFP2} B.~Fantechi, M.~Franciosi, R.~Pardini, {\em Deformation of semi-smooth varieties}, International Mathematics Research Notices {\bf 23}  (2023), 19827?19856.

\bibitem{Fu} W.~Fulton, {\em Intersection Theory}. Springer University Press (1984).


\bibitem{HL} D.~Huybrechts and M.~Lehn, {\em The geometry of moduli spaces of sheaves}. Aspekte der Mathematik (1997) and Cambridge University Press (2010).

\bibitem{L} H.~Lange, {\em Universal Families of Extensions}, J. of Algebra {\bf  83} (1983), 101--112.

\bibitem{Mu}  S.~Mukai, {\em Symplectic structure of the moduli space of sheaves on an
abelian or K3 surface}, Invent. Math.{\bf  77} (1984), no. 1, 101?116.

\bibitem{N} A.~Nobile, {\em From formal smoothings to geometric smoothings}, Rend.~Mat.~Appl. {\bf 44} (2023), 181--210.

\bibitem{Tz} N.~Tziolas, {\em Smoothings of schemes with nonisolated singularities}, Michigan Mathematical Journal {\bf 59} (2010), 25?84.

\bibitem{Y} T.~Yonemura, {\em  Hypersurface simple K3 singularities}, T\'ohoku Math. J. {\bf 42} (1990), 351-380.

\end{thebibliography}
\end{document}